\documentclass[10pt,a4paper]{article}
\usepackage[utf8]{inputenc}
\pdfoutput=1
\usepackage{amsmath}
\usepackage{amsfonts}
\usepackage{amssymb}
\usepackage[numbers]{natbib}
\usepackage[a4paper, total={6in, 8in}]{geometry}
\newcommand{\dian}[1]{\ensuremath{\underline{#1}}}

\newcommand{\media}[2]{\ensuremath{\frac{\partial{#1}}{\partial{#2}}}}

\newcommand{\xw}[1]{\ensuremath{{\bf R}^{#1}}}
\newcommand{\xwx}{\ensuremath{{\bf R}}}

\newcommand{\kxwx}{\ensuremath{{\bf C}}}

\newcommand{\kseks}{\begin{eqnarray}}
\newcommand{\tleks}{\end{eqnarray}}
\newcommand{\kseksw}{\begin{eqnarray*}}
\newcommand{\tleksw}{\end{eqnarray*}}

\newcommand{\eps}{\ensuremath{\epsilon}}

\newcommand{\sge}[1]{\ensuremath{\mathbf{S}^{\Gamma}_{#1}(0)}}

\newcommand{\kl}{\ensuremath{\nabla}}

\newcommand{\olwh}[1]{\ensuremath{ \int_{\mathcal{H}\setminus\{\underline{0}\} }}}

\newcommand{\oppro}[1]{\Pi_{\underline{m},\underline{\mu},\underline{\beta}}^{(#1)}}

\newcommand{\upere}{\ensuremath{\mathcal{H}}}
\newcommand{\supere}{\ensuremath{\mathbf{L}}}

\newcommand{\eug}{\ensuremath{\underline{E}_\Gamma}}

\newcommand{\norma}[1]{\left|\left|#1\right|\right|}
\newcommand{\diag}{\mbox{diag}}
\newtheorem{param}{Theorem}[section]
\newtheorem{reas}{Theorem}[section]

\newtheorem{asexp}{Theorem}[section]
\newtheorem{basest}{Lemma}[section]

\author{Demetrios A. Pliakis\\ ADITAL LLC\\{\ttfamily dpliakis AT gmail.com}}
\title{\bf Small time heat expansion of the laplacian on an analytic hypersurface with an isolated singularity}

\begin{document}

\maketitle

\begin{abstract}

I prove the existence of small time heat expansion for the Laplace operator on an analytic hypersurface with an isolated singularity. First we obtain a local parametrization of the hypersurface near the singularity. We introduce the
notion of quasihomogeneous tangent cone. Then  perturb the parametrization of the cone employing a Newton scheme  and obtain a parametrization with functions of specific form. These allow us to obtain local models for the Laplace operator near the singularity. These are operators 
with irregular singularities. We derive the estimates required by  singular asymptotics. 

\end{abstract}
\section{ Introduction}

The small time heat expansion of a positive second order elliptic operator $H$
 defined on a smooth manifold $M$ is used in index calculations of operators that contain
geometric information about the manifold $M$.

Specifically, let $L$ be a first order differential operator of geometric interest, $L^{\dagger}$ its adjoint then the celebrated McKean-Singer
formula states that 
$$\mbox{ind}(L)=\mbox{tr}\left(e^{-tL^\dagger L}-e^{-tLL^\dagger}\right)=
\lim_{t\rightarrow 0^+}\mbox{tr}
\left(e^{-tL^\dagger L}-e^{-tLL^\dagger}\right)$$
The small time expansion of the 
 self adjoint second order differential operators 
operators $L^\dagger L,LL^\dagger $ arise naturally. 
In particular the asymptotic 
coefficients are spectral invariants of $$L^\dagger L,LL^\dagger $$ Cheeger in his resolution of the Ray-Singer conjecture, in the course of topological operations  encountered spaces with {\it conical singularities}, cones over smooth manifols. However the classical heat expansions are no longer valid and are modified with additional logarithmic terms. Callias-Taubes departing from quantum field theory came across with such singular heat expansions. Namely the computed  determinants that arise in the case of fermions in instanton fields.  The main issue is the asymptotics of integrals that are of power-log form. These type of expansions were already  common in algebraic geometry and their existence is asserted through the celebrated Atiyah-Bernstein theorem. The direct appeal to the Atiyah-Bernstein theorem is not always possible, a simple case is achieved in \cite{P}. Therefore direct analytic methods for the treatment of these integrals were developed and are based either 
{\it singular asymptotic lemma} (\cite{C0},\cite{C2},\cite{CU}) or Melrose's \cite{me1} {\it push forward formula}. 

These methods were efficient for singular spaces with cone like singularities. In these cases the operators that appear contain regular singular points and their resolvents are expressed through Bessel functions. 
   Callias obtained the small time expansion 
for operators with irregular singular points \cite{C2}. Here we built on his results and prove the existence of  heat expansion for the laplacian on an analytic  hypersurface with an isolated singularity. Indeed the laplacian on 
a real  analytic hypersurface with an isolated singularity leads to an operator with an irregular singularity. This requires a parametrization of the hypersurface near the singularity. In the case of complex algebraic curves this is given by 
the classical Puiseux expansion and leads again to conical singularities. This was achieved in \cite{bl}. Moreover in the case of varieties defined by quasihomogeneous polynomials there is also such a parametrization, given by the {\it quasihomogeneous  blow up}.  Toric varieties possess also a direct parametrization with monomial maps but is inefficient for our purposes.
In the general  case of an analytic  hypersurface with isolated singularity (abbreviated here as AHIS) it is neccessary to construct a parametrization and we perform this here. This parametrization is a perturbation of the quasihomogeneous blow up. Actually we reduce an AHIS to its quasihomogeneous tangent cones with the reasonable introduction of quasihomeneous scaling. Then 
we introduce appropriate function spaces inspired from the Puiseux series and 
incorporate a suitable Newton method. This construction allow me to calculate the model operator and employ the methods form \cite{C0},\cite{C1},\cite{C2} for the existence of the asymptotic expansion of the distributional trace of the heat operator. The singularity invariants  defined through the Newton diagramm appear in the exponents of the exponents of the expansion.
\section{Notation and results}
Let $\upere\subset \xw{n+1}$ be the germ of an analytic hypersurface at 
the origin, having an isolated  singularity at the origin. In order to avoid
weird situations that are not interesting at the moment we assume  
that the hypersurface is irreducible and Zariski dense in its 
complexification.
Let $f:\xw{n+1}\rightarrow {\bf R} $ be the 
germ at origin of an analytic function that defines $\upere$ for $\delta>0$:
\kseks
\upere\cap \mathbf{B}_\delta(0)=\{ |\dian{x}|<\delta: f(x)=0\} 
\tleks
Then  the {\it Newton polytope}
\kseks
\mathbf{N}(f)=\{ \alpha\in {\bf N}^{n+1}/ f_\alpha=\media{^\alpha f}{x^\alpha}(0)\neq0\}
+{\bf N}^{n+1}
\tleks
with Newton diagramm $\Delta(f)$ i.e. 
the union of the compact faces of ${\bf N}(f)$. If $\Gamma_j $ is a compact 
face of  \(\mathbf{N}(f)\) then 
\[
\Delta(f)=\bigcup_{j=1}^N\Gamma_j
\]
The face $\Gamma$ defines a quasihomogeneous polynomial of type 
$(\sigma_1,\dots,\sigma_{n+1},m_\Gamma)\in {\bf Q}^{n+2}$ and call 
\(\sigma_1,\dots,\sigma_{n+1}\) {\it quasihomogeneity exponents}
\kseks
f_\Gamma(x)=\sum_{\alpha\in \Gamma}f_\alpha x^\alpha
\tleks
since 
\[
\gamma=(\gamma_1,\dots,\gamma_{n+1})\in \Gamma: \gamma_1\sigma_1+\cdots+\gamma_{n+1}\sigma_{n+1}=m_\Gamma
\]
and
\[
f_\Gamma(\lambda^{\sigma_1}x_1,\dots,\lambda^{\sigma_{n+1}}x_{n+1})=\lambda^{m_\Gamma}f_\Gamma(x_1,\dots,x_{n+1})
\]
Its zero set \(\upere_\Gamma\).  Let also 
\kseks
\eug=\sum_{j=1}^{n+1}\sigma_jx_j\frac{\partial}{\partial x_j} 
\tleks
be the corresponding Euler vector fields, 
\[
\eug(f_\Gamma)=m_\Gamma f_\Gamma
\]
Introduce the function: \(\phi:\xw{n+1}\rightarrow \xw{}\):
\[
\phi_\Gamma(x)=x_1^{\frac{2m}{\sigma_1}}+\cdots+
x_{n+1}^{\frac{2m}{\sigma_1}}
\]
that is obviously quasihomogeneous 
of type $(\sigma_1,\dots,\sigma_{n+1},2m)$ 
 Consider further the sets that we call {\it real Brieskorn spheres} 
 \(sge(\epsilon)\) 
\kseks
\sge{\delta} & = & \{ x\in \xw{n+1}: \phi_\Gamma(x):=\delta^{2m}\}
\tleks 
Write also
\kseks
f(x) & = & f_\Gamma(x)+ R_\gamma(x)\\
\tleks
 We will see that the hypersurface 
$\upere$ decomposes  near the origin into semianalytic pieces 
$\supere_\gamma $ that reveal the approximation by the {\it 
quasihomogeneous tangent cones} that we define and provide the local parametrization of the hypersurface. Recall that the degree of the Gauss map of the link  is the {\it Milnor number} \(\mu\) of the singularity and gives a bound for  the number of branches  of the hypersurface that emanate from the origin. 

We introduce also the space of functions \(P(\epsilon,\Omega)\) generalizing Puiseux seriesa nd we specify later. Denote also \(C_\epsilon(\Omega)=[0,\epsilon]\times \Omega\)

\begin{param} 

 Let $\varepsilon_j>0, j=1,\dots,N$. Let \(\Omega_j\subset \xw{n-1}\) 
 then there is a  collection of maps for  :
\[
\Phi_j :C_{\epsilon_j}(\Omega_j)  \rightarrow \xw{n+1}, \qquad \Phi_j\left(C_{\epsilon_j}(\Omega_j)\right)\subset \upere 
\]
 such that \(\Phi_j\in P(\epsilon,\Omega)^{n+1}\) with a sequence of rational exponents \(\{q_{j,\ell,k}\},j=1,\dots,N,\ell=1,\dots,n+1,k=1,2,\dots\) determined by the Newton digramma. 

\end{param}

These maps are essentially a perturbation of the quasihomogeneous
blow-up.  Their  construction is based on the Newton method. The essence
of this method is that it could be generalized to aritrary codimension and 
complicated singularities; this is under preparation at the moment.

Having obtained this parametrization we have also model operators in hand 
that allow us  to prove the existence of the small time heat expansion
for the Laplace-Beltrami $\Delta_\upere$ on $\upere$ treating its Friedrichs self-adjoint extension. Let 
\[
U=\mathbf{B}_\delta(0)
\]
and spilt this as
\[
U=\bigcup_{j=1}^N U_j
\]
Indeed we prove the following 

\begin{reas}

 Let $\chi_j \equiv 1$ in $U_j$ and 
$\chi\equiv 0 $ outside $U_j$. Then for $\chi=\sum_j\chi_j$ we have 
that as $t\rightarrow 0^+$:
 $$\mbox{tr}(\chi e^{-t\Delta_\upere})\sim t^{-\frac{n}{2}}\sum_{k=0}^\infty 
c_k t^k+\sum_{l=0}^\infty\sum_{i=0}^{S(l)}u_{l,i}[\chi]t^{\alpha_l}\log^it$$
where $\{\alpha_l\}$ is an increasing sequence of rational numbers
determined explicitly by the sequences $\{ \mu_{j,l,k}\}$. The map
$S:{\bf N}\rightarrow {\bf N}$ is also determined by  the Newton diagramm of
the singularity. The $c_k$'s are usual heat coefficients while the $u_{l,i}$'s
are distributions with $\mbox{supp}(u_{l,k})=\{ 0\}$.

\end{reas}
\section{The quasihomogeneous case}

Let \(\upere\) be an AHIS with Neton polytope \(\Delta(f)\). We assume that \(f\) is chosen so that no coordinate plane is contained in \(\upere\). This possible after a suitable rotation. We will consider first the case of 
\(\Delta(f)\) consisting of a single face \(\Gamma\). Then we deal with a {\it quasihomogeneous hypersurface singularity}, \(upere\). Then we have the vector \((\sigma_1,\dots,\sigma_{n+1})\in \mathbf{N}^{n+1}\) and \(m\) such that for all \(\gamma\in \Gamma\) 
\[
\gamma_1\sigma_1+\cdots+\gamma_{n+1}\sigma_{n+1}=m
\]
Then consider the interesection with real Brieskon spheres, \(\epsilon>0\):
\[
\supere(\epsilon)=\sge{\epsilon}\cap \upere
\]
is smooth for all \(\epsilon>0\) due to the fact that on \(\supere(\epsilon)\)
\[
\eug(f)=0 \qquad \eug(\phi_\Gamma)=2m\epsilon^{2m}>0
\]
and hence \(\upere\) and \(\sge(\epsilon)\) meet transversely.

Then write look for \(\xi\) such that 
\[
f(\xi)=0 \qquad \phi_\gamma(\xi)=1
\]
and due to trasversality and the implicit function theorem there exist \(\Omega\subset \xw{n-1}\) and maps such that after a possible renaming of the variables
such that for \(\xi'=(\xi_3,\dots,\xi_{n+1})\)
\[
\Psi_1,\Psi_2: \Omega\rightarrow \xwx,\qquad
\xi_1=\Psi_1(\xi'),\xi_2=\Psi_2(\xi') 
\]
Therefore we have that for a mapping \(\zeta:\Omega\rightarrow \xw{n+1},\Omega\subset \xw{n-1}\) if we denote \(\eta=\xi'\)
\[
x_j=r^{\alpha_j}\zeta_{\eta_j}
\]
Then the parametrization is given if we perform \(r\mapsto r^\lambda\) so that we make \(\min\{\sigma_1,\dots,\sigma_{n+1}\}=1\)
\[
x_1=r^{\sigma_1}\Psi_1(\eta), \qquad x_2=r^{\sigma_1}\Psi_2(\eta)
\]
and for \(i>2\) we have that 
\[
x_i=r^{\sigma_i}\eta_{i-2}
\]
Then the metric indeced on \(\upere\) takes the form
\[
g_\upere= \chi(r,\eta)dr^2+\beta(r,\eta) dr+\gamma_{\supere}(r,\eta)
\]
where
\begin{eqnarray*}
\chi(r,\eta)=\sum_{j=1}^{n+1}\alpha_j^2r^{2(\alpha_j-1)}\zeta_j\\
\beta(r,\eta)=2\sum_{j=1}^{n+1}\sum_{k=1}^{n-1}\alpha_j r\zeta_j\zeta_{j,k}d\eta_k\\
\gamma(r,\eta)=\sum_{i,j=1}^{n-1}\gamma_{k\ell}(r,\xi)d\eta_kd\eta_\ell\\
\gamma_{k\ell}(r,\eta)=r^2\sum_{j=1}^{n+1}r^{\alpha_j+\alpha_k-2}\zeta_{j,k}\zeta_{j,\ell}
\end{eqnarray*}
Notice that  introducing the \((n+1)\times (n-1)\) matrix 
\[
\Lambda_{jk}=\zeta_{j,k}
\]
as well as the \((n+1)\times (n+1)\) diagonal matrix
\[
D_{jj}=\alpha_jr^{\alpha_j-1}, \qquad d\eta=(d\eta_1,\dots,d\eta_{n-1})
\]
we have that 
\begin{eqnarray*}
\chi(r,\eta)= || D\zeta||^2\\
\beta(r,\eta)=2(\zeta, D\Lambda d\eta)\\
\gamma(r,\eta)=r^2\Lambda^TD^2\Lambda
\end{eqnarray*}
\section{Puiseux functions}

Let \(\Omega=[-\delta,\delta]^n\subset \xw{n},\delta>0\) be a cube and \(\epsilon>0\) and form the cylinder \(C_\epsilon(\Omega)=[0,\epsilon]\times \Omega\subset  \xw{n+1}\). Let 
\(\{q_n\}\subset \mathbf{Q}_+\) be an increasing   sequence of nonnegtive rational numbers  and 
\[f_n:\Omega\rightarrow \xw{}, \Omega \ni \xi\mapsto f(\xi)\]  is a sequence of analytic function such that 
\[
\sum_{n=1} \sup_{\Omega} ||f_n||_{\Omega} \epsilon^{q_n}<\infty 
\]
where an analytic function in \(\Omega\) 
\[
f(\xi)=\sum_{\alpha\in \mathbf{N}^n} f_\alpha x^\alpha
\]
the norm is 
\[
||f||_{\Omega}=\sum_{\alpha\in \mathbf{N}^n} f_\alpha \eta^{|\alpha|}
\]
We call the function 
\[
f(r,\xi)=\sum_{n=1}^\infty f_n(\xi) r^{q_n} 
\]
a {\it Puiseux function} and denote their space \(P(\Omega,\epsilon)\). It is elementary to check that  the sum and products of Puiseux functions are also Puiseux functions by merging the sequences of exponents. We alosn define the norm in \(P(\Omega,\epsilon)\) as
\[
|f|_{P(\omega,\epsilon)}=\sum_{n=1} ||f_n||_{\Omega} \epsilon^{q_n}
\]
and \(P(\Omega,\epsilon)\) is a Banach space. 
Let \(\Omega'\subset \xw{k},\Omega\subset \xw{n}\) domains and  \(g:\Omega'\rightarrow \xw{}\) be a real analytic 
function and \(f_1,\dots,f_k\in P(\Omega,\epsilon)\). Then for 
\(F=(f_1,\dots,f_k)\) the function \(g\circ F:C_\epsilon(\Omega)\rightarrow\xw{}\) is in 
\(P(\Omega,\epsilon)\). Indeed let the partial sums
\[
g_N(x)=\sum_{i=0}^N\sum_{|\alpha|\leq k} g_\alpha x^\alpha
\]
that converges in the polydisc \(Q(\eta)=[-\eta_1,\eta_1]\times [-\eta_k,\eta_k]\). We have also 
\[
x_j^{N_j}=\sum_{s=1}^{N_j} f_{m^j_s}(\xi)r^{q^j_s}
\]
Moreover assume that after possible shrinking \(\Omega'\):
\[
\sup_{\xi\in \Omega}|f_{m^j_s}|\epsilon^{q^j_s}\leq \eta_j \qquad \Omega'\subset 
Q(\eta)
\]
Set also \(s=(s_1,\dots,s_k)\) and
\[
Q: \mathbf{N}^k\times \mathbf{Q}_+^k\rightarrow \mathbf{Q}
\]
\[
Q(\alpha,s)=\alpha_1q^1_{s_1}+\cdots+a_kq^k_{s_k}
\]
We form  an icreasing sequence that we denote as \(\{Q_m\}\) and 
set as \(Q=\max\{q(\alpha,s), s_1\leq N_1,\cdots, s_k\leq N_k\}, Q_{N'}=Q\).
Furthermore 
\[
S_r=\{(\alpha,s): q(\alpha,s)=Q_r  \}
\]
 Then we compute 
\[
g_N(r,\xi)=g_N(r,\xi)=\sum_{r=1}^{N'} G_r(\xi) r^{Q_r}
\]
where
\[
G_r(\xi)=\sum_{(\alpha,s)\in S_r}\prod_{j=1}^k\left(\sum_{\beta^j_1+\cdots+\beta^j_{s_{N_j}}=\alpha_j}g_\alpha
\binom{\alpha_1}{\beta^j_1\cdots \beta^j_{s_{N_j}}}f_{m^j_1}^{\beta^j_{s_{N_j}}}\cdots f_{m^j_{s_{N_j}}}^{\beta^j_{s_{N_j}}}\right)
\]
Hence
\[
\sup_{\xi\in Omega}|G(\xi)|\leq 
\sum_{(\alpha,s)\in S_r}|g_\alpha|\prod_{j=1}^k\left(\sum_{\beta^j_1+\cdots+\beta^j_{s_{N_j}}=\alpha_j}
\binom{\alpha_1}{\beta^j_1\cdots \beta^j_{s_{N_j}}}\sup_{\xi\in \Omega}|f_{m^j_1}|)^{\beta^j_{s_{N_j}}}\cdots \sup_{\xi\in \Omega}(|f_{m^j_{s_{N_j}}}|)^{\beta^j_{s_{N_j}}}\right)
\]
Then 
\[
\sum_{r=1}^\infty \sup_{\xi\in \Omega}|G(\xi)|\epsilon^{Q_r}\leq 
\sum_{\ell=0}^\infty\sum_{|\alpha|=\ell} |g_\alpha||\eta|^\alpha<\infty
\]
\section{Quasihomogeneous tangent cones}

For a face of the Newton diagramm of the singularity \(\Gamma\) with  vector \(\sigma=(\sigma_1,\dots,\sigma_{n+1})\):
\[
\sigma_1\gamma_1+\cdots+\sigma_{n+1}\gamma_{n+1}=m_\Gamma 
\]
for all \(\gamma=(\gamma_1,\dots,\gamma_{n+1})\in\Gamma\). Then we introduce the
scaling operator for \(t>0\):
\[
S_{t,\sigma}(x)=\left(t^{\sigma_1}x_1,\dots,t^{\sigma_{n+1}}x_{n+1}\right)
\]
Then we introduce as {\it quasihomogeneous tangent cone} \(K_\Gamma(\upere)\) of \(\upere\) at the singular point at \(0\) as:
\[
K_\Gamma(\upere)=\{\eta\in \xw{n+1}: \eta=\lim_{S_{t_j,\alpha}\rightarrow+\infty} t_jx_j, x_j\in \upere\}
\]
We show that \(K_\Gamma(\upere)\) is analytic set: for  \(\eta\in K_\Gamma(\upere)\) then
\[
f_\Gamma(\eta)=0
\]
Write
\[
f(x)=f_\Gamma(x)+R_\Gamma(x)
\]
Then for the sequance \(\{x_j\}\subset \upere\) we have that
\[
f(x_j)=0\Rightarrow f_\Gamma(x_j)+R_\Gamma(x_j)=0
\]
Then since as \(t_j\rightarrow \infty\):
\[
x_j=\eta_j+\rho_j, \eta_j=S_{1/t_j,\alpha}\eta, \rho_j=S_{1/t_j,\alpha}\rho,\qquad |\rho|=o(1)
\]
we have that due to the analyticity of \(f\) and Taylor series
\[
\frac{1}{t_j^{m_\Gamma}}\left(f_\Gamma(\eta)+\rho\right)+
\frac{1}{t_j^{m_\Gamma'}}\left(R_\Gamma(\eta)+\rho\right)=0\Rightarrow
f_\Gamma(\eta)+o(1)+
\frac{1}{t_j^{m_\Gamma'-m_\Gamma}}\left(R_\Gamma(\eta)+o(1)\right)=0
\]
Since by convexity of the Newton diagramm \(m_\Gamma'>m_\Gamma\) taking the 
limit as
\[
f_\Gamma(\eta)=0
\]
and consequently as well 
\[
f_{\Gamma}(\eta_j)=0
\]
Also we have that for \(\eta\in K_\Gamma(\upere), x_k=\eta_k+\rho_k\) 
\[
 |f(\eta)|\leq \frac{c_1}{t_k^{m'-m}} \qquad  
 |f(\eta_k)|\leq \frac{c_1}{t_k^{m'}}
\]
and selecting \(\epsilon>0\)
\[
\frac{c_1}{t_k^{m'-m}}<\epsilon
\]
Also we have that
\[
|f(\eta)|\leq c_1\epsilon\qquad |f(\eta_k)|\leq c_1\epsilon
\] 
Similarly
\[
\nabla f(\eta)=\nabla f_\Gamma(\eta)+v(\eta)\qquad (1)
\]
for \(c_2>0\):
\[
|v(\eta)|\leq c_2\epsilon
\]
Therefore denoting \(g_j=\frac{\partial g}{\partial x_j}\):
\[
<\nabla f(x_k),\nabla f_\Gamma(x_k)>^2=|\nabla f(x_k)|^2|\nabla f_\Gamma(x_k)|^2-\sum_{i<j}\left(f_i(\eta)f_{\Gamma,j}(x_k)-
f_j(x_k)f_{\Gamma,i}(x_k)\right)^2
\]
Now for 
\[
n_{\Gamma,i}=\frac{f_{\Gamma,i}}{|\nabla f_\Gamma|},\qquad
r_{\Gamma,i}=\frac{R_{\Gamma,i}}{|\nabla f|}
\]
 we have through \((1)\):
\[
\sum_{i<j}\left(f_i(x_k)f_{\Gamma,j}(x_k)-
f_j(x_k)f_{\Gamma,i}(x_k)\right)^2=|\nabla f_\Gamma|^2|\nabla f|^2
\sum_{i<j}\left(r_i(x_k)n_{\Gamma,j}(x_k)-
r_j(x_k)n_{\Gamma,i}(x_k)\right)^2
\]
Since
\[
|r_{\Gamma,j}(x_k)|\leq \frac{c}{t_k^{m'-m}} 
\]
that 
\[
|<\nabla f(x_k),\nabla f_\Gamma(x_k)>|\geq (1-\epsilon^2)|\nabla f(x_k)||
\nabla f_\Gamma(x_k)|
\]
Continuity of \(f,f_\Gamma\)  allows us to conclude that for \(\epsilon>0\) there is \(\delta\) such that for \(x\in \xw{n+1}, |f_\Gamma(x)|\leq \delta\) that 
\[
|<\nabla f(x),\nabla f_\Gamma(x)>|\geq (1-\epsilon^2)|\nabla f(x)||
\nabla f_\Gamma(x)|
\] 
Notice that the analytic set \(K_\Gamma(\upere)\) does not have neccessarily  isolated singularities! However the semianalytic set defined by 
\[
K_\Gamma^{\epsilon,\delta}(\upere) =
K_\Gamma(\upere) \cap \{ x\in \xw{n+1}: |f(x)|<\delta\}\cap \{x\in \xw{n+1}: \phi_\Gamma(x)\leq epsilon^{2m}
\]
has an isolated singularity at the origin. We refer to this set as the {\it quasihomogeneous tanget cone} of height \(\delta\), \(d_\eps\)- close to \(\upere\), where \(d_\eps=O(\eps^{m_\Gamma}\). Similarly we have that 
\[
\upere_\Gamma^{\epsilon,\delta}=
\upere \cap \{ x\in \xw{n+1}: |f_\Gamma(x)|<\delta\}\cap \{x\in \xw{n+1}: \phi_\Gamma(x)\leq epsilon^{2m}
\]

\section{The Newton scheme}

We set up a Newton scheme that provides the parametrization of the analytic hypersurface near its singular point. Specifically we will obtain a perturbation of the parametrization of the quasihomeneous singularity. The perturbation lies in \(P(\epsilon,\Omega)\) for suitable \(\epsilon,\Omega\). 
Let \(\eta\in K_\Gamma^{\epsilon,\delta}(\upere)\) then we set 
\[
\eta=S_{r,\alpha}(\zeta),\qquad \zeta \in \supere_\Gamma^{\epsilon,\delta}=K_\Gamma^{\epsilon,\delta}(\upere)\cap \sge{\epsilon}
\]
which is a smooth analytic hypersurface with an analytic parametrization for 
\(\Omega\subset \xw{n-1}\): 
\[
\zeta: \Omega\rightarrow \xw{n+1}: \xi\mapsto \zeta(\xi),\qquad f(\zeta(\xi))=
\phi_\Gamma(\zeta(\xi))=\delta^{2m}
\]
Then we introduce the map
\[
N(t,\eta)=t-\frac{f(\eta+tn_\Gamma(\eta))}{f'(\eta+tn_\Gamma(\eta))}
\]
for
\[
n_\Gamma(\eta)=\frac{\nabla f_\Gamma(\eta)}{|\nabla f_\Gamma(\eta)}
\]
Notice that since
\[
f_{\Gamma,j}(S_{r,\alpha}\xi)=r^{m-\alpha_j}f_{\Gamma,j}(\xi)
\]
then 
\[
n_\Gamma(\eta)
\]
is a regular function in \(C_\epsilon(\Omega)\). Moreover
the map is well defined in \(P(\epsilon,\Omega)\).
\[
N: P(\epsilon,\Omega)\rightarrow P(\epsilon,\Omega)
\]
We establish that \(N\) is a contraction in \(P(\epsilon,\Omega)\). We assume that uniformly in \(\xi\):
\[
t=O(r^\mu),\qquad \mu=|\alpha|_\infty\]
 We compute 
\[
N'(t)=\frac{f(t)f''(t)}{(f'(t))^2}
\]
where
\[
f'(t)=\nabla f(\eta+tn_\Gamma(\eta))\cdot n_\Gamma(\eta)
\]
\[
f''(t)=\sum_{i,j=1}^{n+1}f_{ij}(\eta+tn_\Gamma(\eta))
n_{\Gamma,i}(\eta)n_{\Gamma,j}(\eta)
\]
Then by Cauchy-Schwarz and the estimate introduced above
\[
|N(t)|\leq \frac{|f(t)||H(f)(\eta+tn_{\Gamma})}{(1-\epsilon^2)|\kl f(\eta+tn_{\Gamma})|^2}
\]
Notice that this again regular since for \(t=O(r^\mu)\)
\[
|H(f)|\leq Cr^{m-2\mu}
\]
Then since for suitable choice of \(\epsilon,\delta\) for the link \(\supere_\Gamma^{\epsilon,\delta}\)
\[
|f(\eta+tn_{\Gamma}(\eta)|\leq r^m\left(|f_\Gamma(\eta)|+r^{m'-m}R_\Gamma(\eta)|\right)\leq 
\delta  r^m
\]
for \(\delta<1\). We look up for function with exponents given by linear combinations of the form for the quasihomogeneity exponents
\(\sigma_1,\dots,\sigma_{n+1}\) 
\[
\alpha_1\sigma_1+\cdots+\alpha_{n+1}\sigma_{n+1} \qquad \alpha_1,\dots,\alpha_{n+1}\in \mathbf{Z}_{\geq 0}
\]

\[
T_n(r,\xi)=\sum_{\ell=1}^\infty T_{n,\ell}(\xi) r^{q_n(\ell)} 
\]
then 
\[
T_{n+1}(r,\xi)-T_n(r,\xi)=\phi'(T_n)(T_n(r,\xi)-T_{n-1}(r,\xi))+
\frac12\phi''(\xi_n)(T_n-T_{n-1})^2
\]
Then setting \(\zeta_{n}(r,\xi)=T_n(r,\xi)-T_{n-1}(r,\xi)\) we have that 
\[
\zeta_{n+1}=\phi'(T_{n-1})\zeta_n+\frac12\phi''(\xi_n)\zeta_n^2 \qquad \mbox{IR}
\]
Substituting and collecting terms we obtain for \(c<1\) and the uniform norm:
\[
||\zeta_{n+1,1}||\leq c||\zeta_{n,1}||
\]
and for \(\ell>1\):
\[
||\zeta_{n+1,\ell}||\leq c||\zeta_{n,\ell}||+F(||\zeta_{n,1}||,\dots,||z_{n,\ell-1}||)
\]
This leads to  recursive inequalities of the form
 \[
 ||\zeta_n||\leq c||\zeta_{n-1}||+\kappa_1 nc^{n}
 \]
 that leads to
 \[
 ||\zeta_n||\leq \kappa_2 n^2c^n
 \]
 and hence that for \(n>m\):
 \[
 ||T_n-T_m||\leq \kappa_2\sum_{k=m}^\infty kc^k\rightarrow 0
 \]
Furthemore  we differentiate (IR)  and obtain for any \(\alpha\in \mathbf{N}^{n-1}\)
\[
D^\alpha \zeta_{n+1}=\left(\phi'(T_{n-1}\frac12\phi''(\xi_n)\zeta_n^2\right)
D^\alpha\zeta_n+F\left(\zeta_n,D\zeta_n,\dots,D^{\alpha-1}\zeta_n\right)
\]
This in turn leads to the following sereis of inequalites for the uniform norm
\[
||D^\alpha \zeta_{n+1}||\leq c_2
||D^\alpha\zeta_n||+F\left(||\zeta_n||,||D\zeta_n||,\dots,||D^{\alpha-1}\zeta_n||\right)
\]
and get for \(k\) depending on \(\alpha\) 
\[
||D^\alpha \zeta_{n+1}||\leq c^n
||D^\alpha\zeta_n||+\kappa_\alpha c^n n^k 
\]
Thne we get the desired convergence of \(D^\alpha T_n\) and therefore thta th limit belongs indeed in \(P(\epsilon,\Omega)\).
  Therefore we have that \(T_n\) convegres in \(P(\epsilon,\Omega)\).
 
 We have shown that if  
\(
 x\in \upere_\Gamma^{\epsilon,\delta}\) there exist  \(\xi,T(r,\xi)\) such that 
\[
x=\chi(r,\xi), \qquad \chi(r,\xi)=\zeta(\xi)+T(r,S_{r,\alpha}\zeta(\xi)))n_\Gamma(\zeta(\xi))
\]
Then  the parametrization is derived setting
\[
\alpha_i=\min\{\alpha_1,\dots,\alpha_{n+1}\}
\]
and
\[
\nu_j=\begin{cases}
 \frac{\alpha_j}{\alpha_i},\qquad \mbox{for} j\neq i\\
 1,\qquad \mbox{for} j=i
\end{cases}
\]
and 
\[
x_k= r^{\nu_k}\chi_k(r,\eta)
\]
where \(\chi_k\in P(\epsilon,\Omega)\). Moreover we have that for \(f_{\Gamma,j}\neq 0\) we have that for \(m\neq j\)
\[
\chi_m(0,\eta)=\eta_m
\]
\section{The Laplacian operator}

\subsection{The metric model}
Given the parametrization obtained above  
\[
C_\epsilon(\Omega) \rightarrow \upere
\]
that has the form with \(\zeta_j\in P(\epsilon,\Omega)\):
\[
x_j=r^{\nu_j}\zeta_j(r,\eta)\qquad j=1,\dots,n+1
\]
and uniformly in \(\eta\)  as \(r\rightarrow 0^+\):
\[
 \zeta_j(r,\eta)=\hat{\zeta}_j+r_j(r,\eta)\qquad r_j(r,\eta)=o(1) 
\]
we compute
\[
dx_j=\nu_j r^{\nu_j-1}\chi_j dr+\sum_{k=1}^{n-1}r^{\nu_k}\zeta_{j,k}d\eta_k,
\]
where
\[
\chi_j(r,\eta)=\zeta_j(r,\eta)+\frac{1}{\nu_j}r\hat{\zeta}_{j,r}
\]
The induced metric on the hypersurface is then
\[
g=\omega dr^2+\sum_{i=1}^{n-1}\beta_jd\eta_jdr+\Sigma
\]
where
\begin{eqnarray*}
A=\diag(\nu_1,\dots,\nu_{n+1})\\
R=\diag(r^{\nu_1-1},\dots,r^{\nu_{n+1}-1})\\
\zeta=(\zeta_1,\dots,\zeta_{n+1})\\
\chi=(\chi_1,\dots,\chi_{n+1})\\
\beta=(\beta_1,\dots,\beta_{n-1})\\
\Lambda_{k,j}=\zeta_{k,j}, k=1,\dots,n+1,j=1,\dots,n-1\\
\omega(r,\eta)=||AR\chi||^2\\
\beta(r,\eta)=r\Lambda(r,\eta)^TAR^2\chi\\
\Sigma(r,\eta)=r^2\Lambda^TR^2\Lambda
\end{eqnarray*}
Notice that \(\lambda>0\) for \(r>0\). We will modify further the metric removing the cross term. 
Therefore solving the system in \(\mathcal{P}\) asking for  the flow out of the link through the initial condtions
\[
\eta_i(\epsilon)=\theta_i
\]
as
\[
\eta_j=\phi_j(r,\theta)
\]
This leads to the singular non-linear system:
\[
\eta'=-\Sigma^{-1}R\beta, \qquad \eta'=\frac{d\eta}{dr}
\]
if we set \(M=R\Lambda, \Sigma=M^TM\). Actually it is an elementary fact that for \(\mu=\max\{\nu_1,\dots,\nu_{n+1}\}-1\) aand since for \(v\in \xw{n-1}\) 
\[
||\Sigma v||_2^2=||R\Lambda v||_2^2\leq cr^{2\mu}||v||_2^2
\]
\[
||\Sigma^{-1}\theta||_q\leq Cr^{-2\mu}||\theta||_p
\]
where \(p,q\geq 1, q=\frac{p}{p-1}\).

This system has a unique solution for \(r>0\) and we will show that the solution extnds to \(r=0\). 
Introduce the function
\[
\phi(\zeta)=\zeta_1^{2m/\nu_1}+\dots+\zeta_{n+1}^{2m/\nu_{n+1}}
\]
Notice that for suitable choice of \(\epsilon>0\) and \(r\in [0,\epsilon]\) we have that 
\[
\phi(\zeta)> c||\eta||_2^{2m}
\]
Moreover we compute
\[
\phi'=2m\left[\frac{1}{\nu_1}\zeta_1^{2m/\nu_1-1}\zeta_1'+\dots+
\frac{1}{\nu_{n+1}}\zeta_{n+1}^{2m/\nu_{n+1}-1}\zeta_{n+1}'\right]
\]
We have that  
\[
\zeta'=\frac{\partial \zeta}{\partial r}-\Lambda \Sigma^{-1}R\beta
\]
Noting that 
\[
\left|\frac{\partial \zeta}{\partial r}\right|\leq C|\Lambda \Sigma^{-1}R\beta|
\]
we invert time as \(t=\epsilon-r\) and consider the equation
\[
\dot{\zeta}=\frac{\partial \zeta}{\partial t}+\Lambda \Sigma^{-1}R\beta\qquad \dot{\zeta}=\frac{d\zeta}{dt}
\]
Apply H\"older inequality and arrive at 
\[
\left|\frac{1}{\nu_1}\zeta_1^{2m/\nu_1-1}\zeta_1'+\dots+
\frac{1}{\nu_{n+1}}\zeta_{n+1}^{2m/\nu_{n+1}-1}\zeta_{n+1}'\right|\leq
2m\left[\zeta_1^{(2\lambda_1-1)p}+\dots+
\zeta_{n+1}^{(2\lambda_{n+1}-1)p}\right]^{1/p}||\zeta'||_q
\]
for  
\[
\lambda_j=\frac{m}{\nu_j}
\]
Now for suitable choice of \(\delta>0\) depending on \(\epsilon\) we have that 
\[
||\dot{\zeta}||_p\leq C
||\Lambda \Sigma^{-1}R\beta||\leq \delta(\epsilon-t)^{-2\mu_+2}||R\beta||_p\leq 
C\delta(\epsilon-t)^{-2\mu+4}||\zeta||_p
\]
Then we select for suitable \(e>0\)
\[
p=\frac{2\Lambda}{2\Lambda-1}+e \qquad \Lambda=\frac{m}{\nu},\qquad \nu=\max\{\nu_1,\dots,\nu_{n+1}\}
\]
or
\[
p=1+e+\kappa, \qquad \kappa=\frac{1}{2\Lambda-1} q=1+\frac{1}{\kappa+e}
\]
and implies 
\[
(2\lambda_j-1)p>2\lambda_j+e(2\lambda_j-1)
\]
This in turn for \(\zeta_j^{2\lambda_j}<\phi\)
\[
\left[\zeta_1^{(2\lambda_1-1)p}+\dots+
\zeta_{n+1}^{(2\lambda_{n+1}-1)p}\right]^{1/p}\leq \phi^{\frac{1+e}{p}}
\]
Also we have that since \(\frac{q}{2\lambda_j}>\frac{q}{2m}\) and some constant \(c>0\)
\[
\left[|\zeta_1|^q+\dots+
|\zeta_{n+1}|^{q}\right]^{1/q}=\left[
\zeta_1^{\frac{2\lambda_jq}{2\lambda_1}}+\dots+
\zeta_{n+1}^{\frac{2\lambda_{n+1}q}{2\lambda_{n+1}}}\right]^{1/q}\leq 
c\phi^{1/2m}
\]
Therefore we end up for \(C\) an \(\delta>0\) to be chosen 
\[
\dot{\phi}\leq C\delta (\epsilon-t)^{-2\mu+4}\phi^{1+\tau}
\]
where for suitable \(e\)
\[
\tau=\frac{1}{2m}-\frac{\kappa}{1+\kappa+e}>0
\]
Then we have that for suitable \(C,\alpha>0\) we have that 
\[
\chi(r)=Cr^\alpha  \qquad \chi(t)=\chi(\epsilon-t), \qquad \dot{\chi}=\frac{\kappa}{(\epsilon-t)^{2\mu-2}} \chi^{1+e}
\]
Then employ the elementary estimmate for
\[
\psi(t)=\frac{\phi(t)}{\chi(t)}
\] 
that satisfies the inequality
\[
\dot{\psi}(t)\leq \frac{\kappa}{(\epsilon-t)^{2\mu-2}}\cdot \psi \cdot \left(\psi^e-1\right)
\]
Then if we arrange \(\epsilon\) so that \(\psi(0)<1\) then 
we have that 
\[
\psi(t)<1
\]
and therefore we have the desired inequality
\[
\phi(\zeta)\leq \delta r^{\alpha}
\]
and hence that the limits exists and vanishes.
In conclusion we have that the metric admits the expression
\[
g=\omega(r,\theta) dr^2+\Sigma(r,\theta)
\]
The hypersurface does not contain any coordinate axis and hence 
\[
\omega(r,\theta)>0
\]
for suitable small \(\epsilon\) and \(0\leq r\leq \epsilon\). Therefore 
we consider the metric
\[
\hat{g}=dr^2+\widehat{\Sigma}(r,\theta)
\]
for 
\[
\widehat{\Sigma}=\frac{\Sigma}{\omega}
\]
\subsection{The operator model}

Employing now the metric near the sngularity we obtain the model oparators near each sector defined by the quasihomegeneous cones where the constants \(\kappa,\mu,\alpha\) have no reference to the constants introduced in the preceding sections while \(\alpha\) depends on the sector 
\[
\Delta_g=\frac{1}{\sqrt{\sigma(r,\theta)}}\partial_r(\sqrt{\sigma(r,\theta)}\partial_r +\Delta_{\widehat{\Sigma}}
\]
where
\[
\sigma=\mbox{det}(\widehat{\Sigma})
\]
Then we will use the model operator model 
\[
H_{\alpha}=H_{0,\alpha,\kappa}+Q_{\alpha}
\]
where for \(k\geq \alpha>2\)
\begin{eqnarray*}
H_{0,\alpha}=\partial_r^2+\frac{\Delta_k}{r^{\alpha}}\\
Q=\frac{L}{r^\alpha}+V(r,\theta)\\
V(r,\theta)=O(r^{-2})\\
\Delta_k=\sum_{i=1}^k\partial^2_{\theta_i}\\
\Delta_{n-1}=\sum_{i=1}^{n-1}\partial^2_{\theta_i}\\
L=L_0-\Delta_k+L_r\\
L_0(\phi)=\frac{1}{\sqrt{\sigma}}\sum_{i,j=1}^k\partial_{\theta_i}\left(\sqrt{\sigma}\Sigma_{ij}\partial_{\theta_j}\phi\right)\\
L_r(\phi)=\frac{1}{\sqrt{\sigma}}\sum_{i,j=k+1}^{n-1}\partial_{\theta_i}\left(\sqrt{\sigma}\Sigma_{ij}\partial_{\theta_j}\phi\right)\\
R_{0,\alpha,\lambda}=\left(H_{0,\alpha}-\lambda\right)^{-1}
\end{eqnarray*}

We start with freezing the coefficients of \(L\) in the \(\theta\) variables  constructing a partition of unity \(\{\chi_\ell\}_{\ell=1}^N\) subordinate to the cver \(\{U_\ell\}_{\ell=1}^N\):
\[
\Omega=\bigcup_{\ell=1}^N U_\ell
\]
such that for \(\Delta_k=\sum_{i=1}^k\partial_{\theta_k}^2\)
\[
||L_0-\Delta_k||\leq \frac{c}{|\lambda|}
\]
This is obtained  using standard estimates for the metric in terms of geodesic normal coordinated in the \(\theta\)-variables and uniform in \(r\). We denote 
\[
H_\alpha=-\partial_r^2-\frac{\Delta}{r^\alpha}
\]
\section{The heat expansion}

\subsection{Operator domain and estimates}
The operator domain \(\mathcal{D}_{0,\alpha}\subset (L^2(C_\epsilon(\Omega))\) that is specified through the following inequality: 

\begin{basest}
 There exist constants \(c,B>0\)  such that for all
\(\phi\in C_0^\infty(C_\epsilon(\Omega))\):
\[
||H_\alpha\phi||^2 \geq ||\partial_r^2\phi||^2+c||\Delta_k r^{-\alpha}\phi||^2-B||\phi||^2
\]
\end{basest}

We follow the lines of Lemma A.1 in \cite{C2}
and compute since \([\Delta,\partial_r]=0\):
\[
(H_\alpha)^2=\partial_r^4+2\Delta_k r^{-\alpha/2}(-\partial_r^2)r^{-\alpha/2}-\frac{\Delta \alpha^2}{2r^{\alpha+2}}+
\frac{\Delta^2}{r^{2\alpha}}
\]
Now use the inequalities \(\phi\in C_0^\infty(C_\epsilon(\Omega))\):
\[
\int_\xwx r^{-2-\alpha}\phi^2\leq C(\alpha) \int_\xwx \left(\partial_r\left(r^{-\alpha/2}\phi\right)\right)^2 
\]
\[
\int_\Omega \phi^2\leq C^2\int_\Omega |\kl \phi|^2 \qquad  \int_\Omega |\kl \phi|^2\leq C\int_\Omega |\Delta \phi|^2 
\]
and since \(\alpha>2\) 
\[
r^{-\alpha-2}\leq \varepsilon r^{-2\alpha}+B
\]
and the inequality follows.

This defines the unique self adjoint extension of \(H_\alpha\) as an unbounded operator in \(L^2\left(C_\epsilon(\omega)\right)\) Then along the same lines of Proposition A.2 we have that  for 
\(R_\alpha(\lambda)=(\lambda-H_\alpha)^{-1}\)
\[
B_{\beta,d}(\lambda)= r^{-\beta}\Delta_k\partial_r^dR_\alpha(\lambda)
\]
and we have that 
\[
||B_{\beta,d}(\lambda)||\leq C_\varepsilon |\lambda|^{-1+\frac{\beta}{\alpha}+\frac{d}{2}}
\]
for \(\Re(\lambda)<\varepsilon |\Im\lambda|-\varepsilon\)
and 
\[
\beta\geq 0 \qquad \frac{\beta}{\alpha}+\frac{d}{2}\leq 1
\]
These lead to the estimate
\[
P_{\alpha}=QR_{0,\alpha,\lambda},  \norma{P_\alpha}\leq c|\lambda|^{-\kappa}
\]
where
\[
\kappa=1-\frac{\beta}{\alpha},\qquad \beta=\max\{\nu_{i,j}\}
\]

\subsection{Neumann series}

We have the standard formula for the heat kernel for \(p\) chosen so that 
the operator inside the integral is trace class: 
\[
e^{-tH_\alpha}=\int_Ce^{-t\lambda} t^{-p}\partial_\lambda^p\left[R_\lambda\right]
\]
for \(R_\lambda=R_\lambda(H_\alpha)\). We write   
\[
R_\lambda=R_\alpha(\lambda)\left[1+Q\cdot 
R_\alpha(\lambda))\right]^{-1}
\]
The preceding estimates allow us to use Neumann series and write for 
\[
R_\lambda=R_\alpha(\lambda)\left[1+PR_\alpha(\lambda)\right]^{-1}=
\sum_{j=0}^\infty R_\lambda^{(j)}
\]
where
\[
R_\lambda^{(j)}=
R_\alpha(\lambda)\left(PR_\alpha(\lambda)\right)^j
\]
Hence 
\[
\partial_\lambda^p\left[R_\lambda\right]=
\sum_{j=0}^\infty\partial_\lambda^p\left[R_\lambda^{(j)}]
\right]
\]
Now for \(c(p,m_0,\dots,m_j)\) being the multinomial coefficient appearing in Leibniz formula for the \(p\)-th derivative of a product
\[
\partial_\lambda^p\left[R_\lambda^{(j)}\right] =\sum_{m_0+\cdots+m_j=p+j+1} c(p,m_0,\dots,m_j)
\Pi_{(m_0,\dots,m_j)}
\]
and 
\[
\Pi_{m_0,\dots,m_j}=R_\alpha(\lambda)^{m_0} Q\cdots
R_\alpha(\lambda)^{m_{j-1}} Q R_\alpha(\lambda)^{m_j} 
\]
Consequently \(\Pi_{m_0,\dots,m_j}\) is written as sum of terms of the form that we denote by \(\Pi_\lambda^j\) and following \cite{C1},\cite{C2}  we call them {\it resolvent products}:
\[
R_\alpha(\lambda)^{m_0} \cdot V_1 \cdot R_\alpha(\lambda)^{m_2} \cdots V_j R_\alpha(\lambda)^{m_j}
\]
where
\[
V_j=r^{\mu_j}f_j(\theta)\partial_\theta^{\beta_j}
\]
with \(\mu_j\in \mathbf{Q}\) while \(\beta_j\) is a polyindex with \(|\beta_j|\leq 2\) and \(f_j\) an analytic function in the \(\theta\) variables.  We denote such a term as
\[
\Pi_{\underline{m}^j,\underline{\mu}^j,\underline{\beta}^j}
\]
with 
\[
\underline{m}^j=(m_0,\dots,m_j),\qquad
\underline{\mu}^j=(\mu_1,\dots,\mu_j)\qquad
\underline{m}^j=(\beta_1,\dots,\beta_j) 
\]
and in concse cotataion
\[
\oppro{j}
\]
We denote as
\[
F_j(t,r,\theta)=t^{\frac{n}{2}-p}\int_C \frac{d\lambda}{2\pi i} e^{-t\lambda} \Pi_\lambda^{(j)}(r,\theta)
\]

For each of these integrals we prove the following theorem

\begin{asexp}

Let \[f_j(\xi,\eta,\theta)\] be defined for \(\xi>,\eta>0,\theta\in \Omega\) then
for
\[
F_j(\xi,\eta,\theta)= f_j(\xi\eta)^\alpha,\eta,\theta)
\]
Then \[F_j\in \Gamma^{S_1,S_2}(\xwx_+\times\xwx_+|\Omega)\]  for
\[
S_1(\ell\alpha)=1=S_2\left(\ell\left(\frac{(3(p-1)+2j)\alpha}{2}-1)\right)\right) \qquad \ell=k,  k=0,1,2
\]
while also \(S_1(z)=1\) for \(z\) non-negative integer linear combinations of the \(\mu_1,\dots,\mu_j\).  Furthermore
\[
P^\eta_z[S_2]F_j(\xi,\eta)|\leq C(\xi\eta)^{\Re z}
\]
\end{asexp}

\subsection{Proof of the asmptotic expansion}

We proceed to the asymptotic expansion of the generic term appearing in the Neumann series and for a constant \(q>\) that we will select later
\[
F_j(t,r,\theta)=t^{\frac{n}{2}+2p}\int_C \frac{d\lambda}{2\pi i} e^{-t\lambda}
\Pi_\lambda^{(j)}(r,\theta)
\]
and obtain the estimates required by the Singular asymptotics lemma \cite{C0},\cite{C2} and we recall in the appendix along with the neccessary notions.
Following \cite{C2} we introduce the variables
 \[
\xi=\frac{t^{1/\alpha}}{r}, \qquad \eta=r
\]
and 
\[
\eta\partial_\eta=\alpha t\partial_t+r\partial_r\qquad 
\xi\partial_\xi=\alpha t\partial_t
\]
Then we have for 
\[
F_j(\xi,\eta,\theta)=F\left(\frac{t^{\frac{1}{\alpha}}}{r},r,\theta\right)
\]
Then we apply 
\[
\eta\partial_\eta F_j(\xi,\eta)=\left(\alpha t\partial_t+r\partial_r\right)\left(
F_j\left(\frac{t^{\frac{1}{\alpha}}}{r},r,\theta\right)
\right)=t^{\frac{n}{2}+q-p}\int_C \frac{d\lambda}{2\pi i} e^{-t\lambda} D\Pi_\lambda^{(j)}(r,\theta)
\]
for
\[
D=\alpha\left(\frac{n}{2}+q-p\right)-\alpha(\lambda\partial_\lambda+1)+r\partial_r
\]
Then 
\[
D\Pi_\lambda^{(j)}(r,\theta)=\sum_{k=0}^j\oppro{k}\left[
\alpha\left(
\frac{n}{2}+q-p\right)+
m_k\left[-\alpha \lambda \partial_\lambda R_\alpha(\lambda)+r\partial_rR_\alpha(\lambda)\right]R_\alpha(\lambda)\right]\oppro{j+1-k}+\left(\mu_1+\cdots+\mu_j\right)\Pi_\lambda^{(j)}+E_j
\]
where the error is computed through the commutator
\[
[r,R_\alpha(\lambda)]=2R_\alpha(\lambda)\partial_rR_\alpha(\lambda)
\]
It introduces a fatcor with norm decaying faster  by a factor of \(|\lambda|^{-\frac12}\) as \(\lambda\rightarrow\infty\).
The important identity from \cite{C2} is modified trivially:
 \[
H_{0,\alpha}=-\frac{1}{2}[r\partial_r,H_{0,\alpha}]+\frac{\alpha-2}{2}\frac{\Delta}{r^\alpha}
\]
and
\[
\lambda\partial_\lambda R_\alpha(\lambda)=-R_\alpha(\lambda)-R_\alpha(\lambda) H_{0,\alpha}R_\alpha(\lambda)
\]
These combine to the formula
\[
\frac{\alpha}{2}R_\alpha(\lambda)-\alpha(\lambda\partial_\lambda R_\alpha(\lambda)+R_\alpha(\lambda))+r[\partial_r,R_\alpha(\lambda]=
\]
\[
= \frac{\alpha}{2}R_\alpha(\lambda)+\alpha R_\alpha(\lambda)H_{0,\alpha}R_\alpha(\lambda)+
[r\partial_r,R_\alpha(\lambda)]-R_\alpha(\lambda)]=
\]
\[
=\left(\frac{\alpha}{2}-1\right)\left[R_\alpha(\lambda)-[r\partial_rR_\alpha(\lambda)]-R_\alpha(\lambda)\frac{\alpha\Delta}{r^\alpha}R_\alpha(\lambda)\right]
\]
\[
=\alpha\left(\frac{\alpha}{2}-1\right)
R_\alpha(\lambda)\partial_rR_\alpha(\lambda)\frac{(-\Delta)}{r^{\alpha+1}}R_\alpha(\lambda)
\]
We obtain then that 
\[
m_k\left[-\alpha\left(\lambda\partial_\lambda R_\alpha(\lambda)+R_\alpha(\lambda)\right)+r\partial_rR_\alpha(\lambda)\right]=
\]
\[
m_k\left[\frac{\alpha}{2}R_\alpha(\lambda)-\alpha\left(\lambda\partial_\lambda R_\alpha(\lambda)+R_\alpha(\lambda)\right)+r\partial_rR_\alpha(\lambda)\right]
-\frac{m_k\alpha}{2}R_\alpha(\lambda)=
\]
\[
=m_k\alpha\left(\frac{\alpha}{2}-1\right)
R_\alpha(\lambda)\partial_rR_\alpha(\lambda)\frac{(-\Delta)}{r^{\alpha+1}}R_\alpha(\lambda)
-\frac{m_k\alpha}{2}R_\alpha(\lambda)
\]
Then we have that
\[
D\Pi_\lambda^{(j)}(r,\theta)=\alpha\left(\frac{n}{2}+q-p-\frac{p}{2}\right)
\Pi_\lambda^{(j)}(r,\theta)+\sum_{k=0}^jm_k\oppro{k}\left[
\alpha\left(\frac{\alpha}{2}-1\right)
R_\alpha(\lambda)\partial_rR_\alpha(\lambda)\frac{(-\Delta)}{r^{\alpha+1}}R_\alpha(\lambda)\right]\oppro{j+1-k}+
\]
\[
+\left(\mu_1+\cdots+\mu_j\right)\Pi_\lambda^{(j)}+E_j
\]
At this point we select  \(p=2\ell,q=3\ell-\frac{n}{2}\) and then employing  commutators introducing factors of with strictly smaler norms we have that for 
a function \(\phi\in L^1,\phi=\phi_1\phi_2,\phi_1,\phi_2\in L^2\) we have tha 
\[
||\phi \oppro{k}\left[ R_\alpha(\lambda)\partial_rR_\alpha(\lambda)
\frac{(-\Delta)}{r^{\alpha+1}}R_\alpha(\lambda)\right]\oppro{j+1-k}||_1\leq C
||\phi 
R_\alpha(\lambda)^\ell\partial_rR_\alpha(\lambda)\frac{(-\Delta)}{r^{\alpha+1}}R_\alpha(\lambda)^{\ell-1}V_1\cdots V_j R_\alpha(\lambda)||
\leq 
\]
\[
c||\phi|| |\lambda|^{-\left(2\ell+j-\frac{n+5}{2}\right)+\frac{1}{\alpha}}
\]
This allows the estimate by duality 
\[
||f||_\infty =\sup_{g\in L^1(\mathbf{R})}\left|\int fg\right|
\]
Then we conclude that 
\[
|\eta\partial_\eta f(\xi,\eta)|\leq C |\xi\eta|^z
\]
for \(z=\left(\frac{3(p-1)}{2}+j\right)\alpha-1\).
\paragraph*{Higher \(\eta\partial_\eta\)-derivatives} These are derived using the commutator technique as is given in \cite{C2}

\paragraph*{\(\xi\partial_\xi\)-derivatives} The \(\xi\partial_\xi\) reduces to \(\alpha t\partial_t\) and hence the \(\xi\rightarrow 0\)asymptotics to small time asymptotics of \(F_j\) that are reduced to classical expansions. 
\newpage
\section{Appendix}

\subsection{Resolvent estimates}

\paragraph{Resolvent comparison formula}

Following the lines of [C] we obtain the resolvent comparison formula: the resolvent of the irregular laplacian to the free laplacian. Let
\(\Delta_n \) be the laplacian in eulidean space and \(R_0(\lambda)\) its resolvent. 
Also let \(\Omega\) be a domain in \(\xw{n-1}\)
\[
\Pi: L^2(\xw{n})\rightarrow L^2((0,\infty)\times \Omega)
\]
Define
\[
R_\lambda^0=\Pi \cdot R_0(\lambda) \cdot \Pi
\]
and
\[
Q_{k,\alpha,\lambda}=\frac{\Delta}{r^\alpha}\cdot R_\alpha(\lambda)
\]
Then we have the formula
\[
R_\alpha(\lambda)=R_\lambda^0+R_\lambda^0\cdot Q_{k,\alpha,\lambda}
\]
\paragraph{Resolvent factors}

We encounter terms of the form
\[
P_{\beta,\gamma,d,\alpha}=r^{\beta}\cdot \partial_\theta^\gamma  \cdot\partial_r^d \cdot R_\alpha(\lambda)^m
(\lambda)
\]
and  its {\it index} as in \cite{C2}:
\[
\mbox{ind}(P_{beta,\gamma,d,\alpha}=\begin{cases}
m-\frac{\beta}{\alpha}+\frac{d}{2}, \quad \mbox{if} \quad \beta \leq 0\\
m-\frac{d}{2}, \quad \mbox{if} \quad  \beta>0
\end{cases}
\]
and its degree \(\mbox{deg}_+\) for \(\beta>0\) as
\[
\mbox{deg}_+=\beta-d+1
\]
Then we have the following estimates that are modified versions of the corresponding ones from \cite{C1}:
\begin{itemize}
\item if \(\beta\leq 0\) then
\[
||\phi P_{\beta,\gamma,d}||_k\leq c||\phi||_{2}|\lambda|^{-\mbox{ind}(P_{beta,\gamma,d}+\frac{n}{2k}}
\]
\item if \(r\in (0,\frac{1}{\sqrt{\lambda}})\) then
\[
||\phi P_{\beta,\gamma,d,\alpha}||_k\leq c||\phi||_{2}|\lambda|^{-\mbox{ind}(P_{0,\gamma,d,\alpha}+\frac{n}{2k}-\frac{\beta}{2}}
\]
\end{itemize}

We recall here the following etimates form \cite{C1} for \(0\leq \alpha\leq 1\)
\[
||\phi r^{-\alpha} R_0(\lambda)^p||_k\leq c_k|\lambda|^{-p-\frac{\alpha}{2}+\frac{n}{2k}}
\]
Employing the preceding estaimeta and the comparison formula and operator norm estimate and derive the estimate
\[
||\phi r^{-\alpha} R_{\alpha}(\lambda)^p||_k\leq c_k|\lambda|^{-p-\frac{\alpha}{2}+\frac{n}{2k}}
\]

\subsection{Singular Asymptotics}

We recall here the basic definitions form [C] neccessary for the sigular asymptotic expansions.
Let \(S:\kxwx\rightarrow\mathbf{Z}_+\) with \(\sum_{\Re(z)<k} S(z)<\infty\) for all \(k\in\kxwx\) and its degree \(m=m(S)=\sum_{z\in \kxwx} S(z)\). Then we introduce the operators:
\[
P_z^x[S]=\prod_{\Re(z')\leq \Re(z),z'\neq z} (x\partial_x-z')^{S(z')}
\] 
Similarly we introduce the space of functions that possess {\it power-log} expasnions at \(0^+\). We call this space \(\Gamma^S(0,\infty)\) and a function 
\[
f:(0,\infty)\rightarrow \mathbf{C}
\]
\(f\in \Gamma^S(0,\infty)\) if for  \(m=m(S)\) \(f\in C^m(0,\infty)\) and 
\[
\partial_t^sf(t)=\sum_{\Re(z)<\Re(k)\sum_[0\leq j<S(z')}f_{zj}\partial_t^s\left[t^z\log^jt\right]+O\left(t^{k+\delta_k-s}\right)
\]
Define similarly the asymptotics for functions of \(n-\)  variables as we approach corners. Then we have the following basic facts:
\begin{itemize}
\item If \(m(S)=\infty\) then \(f\in \Gamma^S(0,\infty)\) iff \(f\in C^\infty(0,\infty)\) and for all \(z\in \mathbf{C}\)
\[
P_z^t[S]f(t)=O(t^{z-\epsilon})
\]
for \(\epsilon>0\).
\item Let \(f\in \Gamma^{S_1,S_2}\left(\xw{}_+\times \xw{}_+\right), m(S_1)=m(S_2)=\infty)\) and for all \(k\in\kxwx\) and \(\epsilon>0\)
\[
\left|P^x_k(S_1)f(x,y)\right|\leq (xy)^{\Re(k)-\epsilon}h_{k,\epsilon}(y)\qquad
\int_1^\infty h_{k,\epsilon}(t)\frac{dt}{t}
\]
Moreover let 
\[
F(t)=\int_0^\infty f\left(x,\frac{t}{x}\right)\frac{dx}{x}
\]
Then 
\[
F\in\Gamma^{S_1+S_2}\left(\xwx_+\right)
\]
\end{itemize}
 
\end{document}